\documentclass[11pt]{article}
\usepackage{amsbsy,amsfonts,amsmath,amssymb,amscd,amsthm}
\usepackage{mathrsfs}
\usepackage[a4paper,text={6.1in,8in},centering,includefoot,foot=1in]{geometry}
\usepackage{amsmath,amsthm,amsfonts,amssymb,amscd}
\usepackage{graphics}
\usepackage{color}
\usepackage[colorlinks=true,citecolor=blue]{hyperref}
\usepackage{latexsym}
\usepackage{cite}
\usepackage{url}
\usepackage[dvips]{graphicx}

\small\normalsize
\theoremstyle{plain}
\newtheorem{theorem}{Theorem}

\newtheorem{lemma}{Lemma}

\newtheorem{proposition}{Proposition}

\numberwithin{equation}{section}

\begin{document}

\title{\textbf{Measure theoretic structure of \\ the  sum of  two  affine  Cantor sets}}
\date{}
\author{  Mehdi  Pourbarat\\
\\} \maketitle

\maketitle

\hspace*{-.5cm}\textbf{Abstract}:  
Suppose that  $K$ and $ K'$ are  two affine Cantor sets in $\mathbb R$.    It is  shown that  the sum set $K+K'$ has equal box and Hausdorff dimensions,  where this common dimension is denoted by $s$, we have   $H^s(K+K')<\infty$. It has previously been proven that   $H^s(K +K') > 0$,   for  almost every pair $(K, K')$ satisfying  $HD(K)+HD(K') <\frac{1}{2}$ or  $HD(K)+HD(K') > 1$.  
We   show  that      $H^s(K +K') = 0$,  for  almost every pair $(K, K')$ satisfying the conditions $HD(K) + HD(K') \leq 1$ and  $\tau(K) + \tau(K') + 3\tau(K)\tau(K') \geq 1$.

\section{Introduction} \label{sec1}
 
Affine Cantor sets  have been regarded, on one hand as   a   natural generalization of the 
 homogeneous  Cantor sets   and   on the other hand as  a particular case of regular Cantor sets. An affine   Cantor set can be obtained via an iteration process on the finite union of many discrete   closed intervals  not necessarily same lengths  \cite{Fa}. In one dimensional  dynamical systems theory,  regular Cantor sets  may  appear   as repellers
of expanding  maps, for definitions and background, see \cite{Fa1}.  In the study of  two-dimensional diffeomorphisms, they also may  appear from intersections  of hyperbolic sets  with stable or unstable manifolds of its points   \cite{PT}. 
 The arithmetic   sum  of two regular Cantor sets  sometimes   play a key role  in various settings such as   dynamical systems, number theory,  real analysis,  spectral theory,  harmonic analysis (for each, see a number of references in \cite{P4}) and
the existence of patterns in fractals \cite{Ma,Y}.

In dynamical systems theory,  a  Cantor set $K$  is \textit{regular} or \textit{dynamically defined} if:
   \begin{itemize}
   	\item [(i)] there are disjoint compact intervals $I_{0}, I_{1},\ldots,I_{M}$ such that
   	$K\subset \bigcup_{i=0}^M I_{i}$ and the boundary of each $I_{i}$ is
   	contained in $K$,
   	\item [(ii)] there is a $C^{1+\epsilon}$ expanding map $\psi$
   	defined in a neighbourhood of set $ \bigcup_{i=0}^M I_{i}$ such that
   	$\psi(I_{i})$ is the convex hull of a finite union of some intervals $I_{j}$
   	satisfying:
   	\begin{itemize}
   		\item [$\bullet$]  for each $0\leq i\leq M$ and $n$ sufficiently
   		large, $\psi^{n}(K\cap I_{i})=K,$
   		\item [$\bullet$]$K=\bigcap_{n=0}^{\infty}\psi
   		^{-n}( \bigcup_{i=0}^M I_{i})$.
   	\end{itemize}
   \end{itemize}
   The set $\{I_{0},I_{1},\cdots,I_{M}\}$ is, by definition, a Markov partition for $K$, 
   the set $D:=\bigcup_{i=0}^{M} I_{i}$ is  Markov domain   of $ K$ and $K_n=\psi
   		^{-n}( \bigcup_{i=0}^M I_{i})$ is called  the $n$-th step of the construction  of $K$.  A regular Cantor set $K$ is \textit{affine} if $D\psi$ is constant on every interval $I_i$. It is called \textit{homogeneous}, if  $D\psi$ is constant on the Markov domain $D$.   A  $C^2$- regular Cantor set $K$ is non \textit{essentially affine}
  if the expanding dynamics which defines it cannot be $C^2$ conjugated to a map whose second derivative is identically zero on the Cantor set. 

Suppose that $K$ and $K'$ are two regular   Cantor sets and $\lambda$ is a real number.  Define
 $$K -\lambda  K' :=\big\lbrace x-\lambda y\mid ~x \in K,~y \in K'\big\rbrace.$$
 Note that  $-K'$ is a regular Cantor set too, thus   the  sum set  $K +\lambda  K'$ can be considered as a  difference set \cite{P11}.  
The reason that we are interested  in considering  difference sets in the current form is that for each element $t\in K -\lambda  K'$ and $\lambda\neq0$,  the sliding image $ t +\lambda  K'$  of  $\lambda  K'$ intersects $K$, and  from this point of view investigating   the geometric structure of the set $K -\lambda K'$  is useful in studying the  bifurcation phenomena in  dynamical systems \cite {PT}. 

In  geometric measure theory, regular Cantor sets may be   defined as  the attractors of the  iterated function systems  (IFSs) which  satisfy the open set condition.  Assume that affine Cantor sets  $K$ and $K'$ are the attractors of the  IFSs $\{F_i(x)=r_ix+t_i\}_{i=1}^{M}$ and $\{F_j'(x)=r_j'x+{t'}_j\}_{j=1}^{N}$ with  convex hulls  $[0,~a]$  and  $[0,~b]$, respectively.  Throughout the paper,  abbreviation  $BD(A),~HD(A),~H^s(A)$,  $|A|$ and $A^\circ $ denote  the  Box dimension,  Hausdorff  dimension, $s$-dimensional  Hausdorff measure, Lebesgue  measure and interior  of $A$, respectively  \cite{Fa1}.

 Moreira \cite{M2} proved that if     $K$ and $K'$ are  two   regular Cantor sets, then 
\begin{equation}
HD(K-K')=\min\{HD(K)+HD(K'), 1\}, 
\end{equation}
provided that  one of Cantor sets  is  non essentially affine. The set of all resulting pairs forms a general  subset  (open and dense subset) of  the space of all pairs of $C^2$-regular  Cantor sets. 
He also showed $BD(K-K')=HD(K-K')$. 
Peres   and  Shmerkin \cite{PS}   proved  (1.1) provided that  $\frac{\log |r_i|}{ \log |r_j'|}\notin \mathbb{Q}$ for at least a pair $(i,j)$.  They also  introduced some examples  to show that in cases where   $\frac{\log |r_i|}{ \log |r_j'|}\in \mathbb{Q}$ for all $(i,j)$,  dimension drop may  occur.  Our first result is as follows:
\begin{theorem}\label{A}  For given  affine Cantor sets  $K$ and $K'$,  
   $BD(K-K')=HD(K-K'):=s$  and  $H^s(K-K')<\infty$. Also,   if   $HD(K)+HD(K')\leq 1$ and  $ r_i,r'_j>0$, moreover  $\frac{\log r_i}{ \log r_j'}\notin \mathbb{Q}$   for at least a pair $(i,j)$,  then   there exists a  dense subset $ D\subset \mathbb{R}$  such that  $H^s(K-\lambda K')=0$ for all  $\lambda \in D$. 
\end{theorem}
We \cite{P4} have already shown the result in the context of    homogeneous  Cantor sets.
Looking at the proof of Theorem \ref{A}, we see that the proof holds even under weaker assumptions.  Indeed, if the IFSs do not  satisfy  the open set condition,  the theorem still holds.
Also,   for all regular Cantor sets $K$ and $K'$ satisfying (1.1), we always have the assertions  $BD(K-K')=HD(K-K'):=s$ and  $H^s(K-K')<\infty$.

To state the second result, we  need to  introduce the concept of    Newhouse's thickness.  Suppose that $K$ is a regular Cantor set and   $U$ be a bounded gap  of  $K$.  Also, suppose that $L_U$ and  $R_U$ are the intervals at its left  and its right, respectively, that separate it from the closest larger gaps;

\unitlength1.30mm \linethickness{1pt}
\begin{picture}(111.33,12.27)
\put(13.00,5.00){\makebox(0,0)[lc]{$K:$}}
\put(20.00,5.00){\line(1,0){5.00}}
\put(35.00,5.00){\line(1,0){10.00}}
\put(53.00,5.00){\line(1,0){6.00}}
\put(70.00,5.00){\line(1,0){5.00}}

\put(20.00,3.70){\makebox(0,0)[cb]{$|$}}
\put(25.30,3.70){\makebox(0,0)[cb]{$($}}
\put(34.70,3.70){\makebox(0,0)[cb]{$)$}}
\put(45.30,3.70){\makebox(0,0)[cb]{$($}}
\put(52.70,3.70){\makebox(0,0)[cb]{$)$}}
\put(59.30,3.70){\makebox(0,0)[cb]{$($}}
\put(69.70,3.70){\makebox(0,0)[cb]{$)$}}
\put(75.00,3.70){\makebox(0,0)[cb]{$|$}}

\put(40.00,7.00){\makebox(0,0)[cb]{$L_U$}}
\put(49.00,5.50){\makebox(0,0)[cb]{$U$}}
\put(56.00,7.00){\makebox(0,0)[cb]{$R_U$}}
\end{picture}
Put $\tau_R(U):=~\frac{|R_U|}{|U|}$ and $~\tau_L(U):=~\frac{|L_U|}{|U|}$.
The \emph{right thickness}  $\tau _{R}$ and \emph{left thickness} $\tau_{L}$ are
$$\tau_R(K):=~\inf_{U}\tau_R(U),~~~~~~~~~~~\tau_L(K):=~\inf_{U}\tau_L(U)$$
and  Newhouse's thickness $\tau(K)$  is the minimum of $\tau_R(K)$ and $ \tau_L(K)$.
\begin{theorem}\label{B}   
Suppose that  $K$ and $K'$ are two  affine Cantor sets  satisfying $HD(K)+HD(K')\leq 1$ and  $\tau(K)+\tau(K')+3\tau(K)\tau(K')\geq1$.  If  $\frac{\log |r_i|}{ \log |r_j'|}\notin \mathbb{Q}$ for at least a pair $(i,j)$,  then 
$H^s(K-\lambda K')=0$ for almost every   $\lambda \in \mathbb R$. 
\end{theorem}
  Marstrand’s projection theorem tells  that if  $HD(K)+HD(K') > 1$, then $K-\lambda K'$ has positive Lebesgue measure;     $H^1(K-\lambda K')>0$, for almost every   $\lambda \in \mathbb R$ \cite {PT}.  Also  from  Theorem 2.14 in \cite{I} and  Theorem 2 in  \cite{PS},  if  $HD(K)+HD(K') <\frac{1}{2}$  and $\frac{\log |r_i|}{ \log| r_j'|}\notin \mathbb{Q}$ for at least a pair $(i,j)$  then 
$H^s(K-\lambda K')>0$, for almost every   $\lambda \in \mathbb R$. Nothing is known in the region where   $HD(K)+HD(K') \geq \frac{1}{2}$ and   $\tau(K)+\tau(K')+3\tau(K)\tau(K')<1$.
If  $C_\lambda$  and $C_\mu$  is the middle-$(1-2\lambda)$ and  middle-$(1-2\mu)$ Cantor
sets, respectively, then  the relation  $\tau(C_\lambda)+\tau(C_\mu)+3\tau(C_\lambda)\tau(C_\mu)\geq1$  is equivalent to  $\frac{2\lambda}{1-3\lambda}\frac{2\mu}{1-3\mu}\geq 1$ and in this case, the result has already been proven by Eroglu \cite{I}.  We also generalize Eroglu's result in the same form as he examined it at the end of his paper.

\section{The proof of Theorem \ref{A} } \label{sec2}
There are two possibilities for similarity ratios:
\begin{itemize} 
\item  If $\frac{\log |r_i|}{ \log | r_j'|}\notin \mathbb{Q}$ for at least a $i$ and $j$, then 

\begin{itemize}
\item [$\bullet$]  
 If $HD(K-K')=1$, then  $BD(K-K')=1$ since $HD(A)\leq \underline{\dim}_B (A)\leq \overline{\dim}_B(A)\leq 1$ for each subset $A$ of $\mathbb R$.  Here,   $\underline{\dim}_B (A)$ and $\overline{\dim}_B(A)$ denote the lower and upper box  dimensions of $A$, respectively \cite{Fa1}. Also, $H^1(K-K')<\infty $, since  $H^1$ coincides with Lebesgue measure.
\item  [$\bullet$] If  $HD(K-K')<1$, then  
\begin{eqnarray*}
HD\big(K- K'\big)&=& HD(K)+HD(K')   ~~~~\text{(Peres-Shmerkin \cite{PS} )} 
\\&=& \overline{\dim}_B(K)+\overline{\dim}_B(K') ~~~~\text{(Product formula 7.5 in \cite{Fa1})}\\&\geq&\overline{\dim}_B(K\times K')~~~~\text{(Exercise 6.12 in \cite{Fa1})}\\&\geq&\overline{\dim}_B(K-K'). 
\end{eqnarray*}

\end{itemize}
Consequently, $BD(K-K')=HD(K-K')$. Also, $H^s(K-K')<\infty $, since 
$$H^{d+d'}(K-K')\leq H^{d+d'}(K\times K')<\infty,$$
where   $d:=HD(K)$ and $d':=HD(K')$,  see Propositions 2.2 in \cite{Fa1}  and 4.4 in \cite{PT}.
\end{itemize}
\begin{itemize} 
\item  If $\frac{\log |r_i|}{ \log |r_j'|}\in \mathbb{Q}$  for all $i$ and $j$, then
\end{itemize}
there are natural numbers $\{m_i\}_{i=1}^{M},~\{n_j\}_{j=1}^{N}$ and  $\gamma<1$ such that  $|r_i|=\gamma^{m_i}$ and $|r'_j|=\gamma^{n_j}$. Let $m=[m_1,m_2,\dots,m_M]$ and  $n=[n_1,n_2,\dots,n_{N}]$, which $[,]$ means  least common multiple (lcm) of the numbers.  We  suppose that $(m_1,m_2,\dots,m_M,n_1,n_2,\dots,n_{N})=1$ which $(,)$ means   greatest common divisor (gcd)  of the numbers. 

\textbf{Claim:} There are  two  finite sequence    $\{\overline{m_i}\}_{i=1}^{M}$ and $\{\overline{n_j}\}_{j=1}^{N}$  of the set ${\mathbb N}_*:=\mathbb N\cup\{0\}$ such that  $\sum_{i=1}^M\overline{m_i}m_i-\sum_{j=1}^N\overline{n_j}n_j=1$. 

For a proof, suppose that integer numbers  $\{\tilde{m_i}\}_{i=1}^{M}$ and $\{\tilde{n_j}\}_{j=1}^{N}$
are the solution of the Diophantine equation $\sum_{i=1}^M\tilde{m_i}m_i-\sum_{j=1}^N\tilde{n_j}n_j=1$. 
If  $\tilde{m_i}\geq 0$, let  $\overline{m_i} =\tilde{m_i}$.  If    $\tilde{m_i}<0$, then there  are  integers  $x_i>0$ and $k_i\geq 0$ such that
$-\tilde{m_i}=x_in_1-k_i$. Thus $\tilde{m_i}m_i=-x_im_in_1+k_im_i$.  Let  $\overline{m_i} =k_i$ and add  $x_im_i$ to $\tilde{n_1}$. Hence, we can suppose that $\overline{m_i}\geq 0$, for all  $1\leq i\leq M$. Similarly, If  $\tilde{n_j}\geq 0$, let  $\overline{n_j} =\tilde{n_j}$.  If    $\tilde{n_j}<0$, then there are  integers $y_j>0$ and $l_j\geq 0$ such that
$-\tilde{n_j}=y_jm_1-l_j$. Thus $-\tilde{n_j}n_j=y_jn_jm_1-l_jn_j$. Let  $\overline{n_j} =l_j$ and add  $y_jn_j$ to $\tilde{m_1}$. Hence, we also  can suppose that $\overline{n_j}\geq 0$, for all  $1\leq j\leq N$ and this completes the proof of the claim

If there are   the compositions of the maps of   $\{F_i(x)=r_ix+t_i\}_{i=1}^{M}$ and $\{F_j'(x)=r_j'x+{t'}_j\}_{j=1}^{N}$ named  $F_K$ and $F_{K'}$, respectively such that $D F_K=-DF_{K'}$, where $D$ is differential operator, then let $c=|DF_K|$, otherwise let $c=1$.
Let  $A:=\frac{c\gamma^{mn+m^*+(m^*+n^*)\sum_{i=1}^M \overline{m_i}m_i}}{a+b}$, where   $m^*=\max\{m_i\}_{i=1}^{M}, ~ n^*=\max\{n_j\}_{j=1}^{N}$  and $R_0:=a+b$.

 For given $t\in K-K'$ and $0<R<R_0$, we introduce a mapping $F:K-K'\to (K-K')\cap (t-R,~t+R)$ satisfying
$$AR|x-y|\leq |F(x)-F(y)|, ~~~~\forall x,y\in K-K'.$$   According  to Theorem 3.2 in  \cite{Fa}, $BD(K-K')=HD(K-K')$ and  $H^s(K-K')<\infty $.

 For each $k\in \mathbb {N}$, let  the sets $M_k:=\{(M_{i,k})_{i=1}^{M}\in \mathbb N_*^M\}$  and $N_k:=\{(N_{j,k})_{j=1}^{N}\in \mathbb N_*^{N}\}$ that  satisfy 
$$kmn-m^*< \sum_{i=1}^{M}m_iM_{i,k}\leq kmn,~~~~~~~kmn\leq \sum_{j=1}^{N}n_jN_{j,k}< kmn +n^*.$$ 
Thus for all  $(M_{i,k})_{i=1}^{M}\in M_k$ and $(N_{j,k})_{j=1}^{N}\in N_k$, we always  have
$$\gamma^{m^*+n^*}<\frac{\Pi_{j=1}^{N}|{r'}_j|^{N_{j,k}}}{\Pi_{i=1}^{M}|r_i|^{M_{i,k}}}=\frac{\gamma^{ \sum_{j=1}^{N}n_jN_{j,{k}}}}{\gamma^{\sum_{i=1}^{M}m_iM_{i,{k}}}}\leq 1.$$
Let $I_{k} $ denote  the set of all  intervals  at the $k$-th step of the construction of $ K$ and  $J_{k} $ denote  the set of all the intervals  at  the $k$-th step of the construction of $ K'$. Let ${\cal{I}}:=\cup_{k\in\mathbb N}I_k$ and  ${\cal{J}}:=\cup_{k\in\mathbb N}J_k$. Finally, let
$$E_k:=\big \{I\times J\in{\cal{I}}\times {\cal{J}}\big | ~kmn-m^*<\log_\gamma\frac{ |I|}{a}\leq kmn,~kmn\leq \log_\gamma\frac{ |J|}{b}< kmn+n^*\big \}.$$
Thus, for each $I\times J\in E_k$, there are  $(M_{i,k})_{i=1}^{M}\in M_k$ and   $(N_{j,k})_{j=1}^{N}\in N_k$ such that 
$|I|=a\Pi_{i=1}^{M}|r_i|^{M_{i,k}}$ and $|J|=b\Pi_{j=1}^{N}|{r'_j}|^{N_{j,k}}$. Moreover,
$$|\pi_{\frac{\pi}{4}}(I\times J)|=a\Pi_{i=1}^{M}|r_i|^{M_{i,k}}+b\Pi_{j=1}^{N}|{r'_j}|^{N_{j,k}}=\Pi_{i=1}^{M}|r_i|^{M_{i,k}}(a+\frac{\Pi_{j=1}^{N}|{r'}_j|^{N_{j,k}}}{\Pi_{i=1}^{M}|r_i|^{M_{i,k}}}b).$$
where  $\pi _{\frac{\pi}{4}}(x,y)=x-y$, (indeed, the  map $\pi _{\frac{\pi}{4}}$ projects the points of  $\mathbb{R}^2$   onto $\mathbb{R}\times \{0\}$ in the  direction of lines with  angle $\frac{\pi}{4}$). It is not hard to show that 
\begin{itemize} 
\item for each $k\in\mathbb N$ and $x\in K-K'$, there is a  $I\times J\in E_k$ such that $x\in\pi_{\frac{\pi}{4}}(I\times J)$,
\item for each $k\in\mathbb N$ and $I\times J\in E_k$,  there are $I_1\times J_1\in E_{k-1}$ and $I_2\times J_2\in E_{k+1}$  such that  $I_2\times J_2\subset I\times J\subset I_1\times J_1$,
\item $\lim_{k\to \infty} |I_k|+|J_k|=0$,  for each the sequence $\{I_k\times J_k\}_{I_k\times J_k\in E_k}$.
\end{itemize}
Since $t\in K-K'$ and $R\leq R_0$, one can find a   natural number $k$ which is  smallest in the sense that there is a $I\times J\in E_{k+1}$ such that $\pi_{\frac{\pi}{4}}(I\times J)\subset (t-R,~t+R)$ and $t\in\pi_{\frac{\pi}{4}}(I\times J)$.  Let $F_I$ and $F_J$  be the compositions of the maps of the IFSs  $\{F_i\}_{i=1}^{M}$ and $\{F_j\}_{j=1}^{N}$, respectively such that $F_I\big([0,~a]\big) =I$ and  $F_J\big([0,~b]\big) =J$. Thus, there is a   $(M_{i,k+1})_{i=1}^{M}\in M_{k+1}$  such that $|I|=a\Pi_{i=1}^{M}|r_i|^{M_{i,k+1}}$. Take a $I_1\times J_1\in E_{k}$  such that  $I\times J\subset I_1\times J_1 $. 
 Thus, there is  a $(M_{i,k})_{i=1}^{M}\in M_{k}$  such that $R< (a+b)\gamma^{\sum_{i=1}^{M}m_iM_{i,k}}$. 
Let $ l:=\log_\gamma\frac{\Pi_{j=1}^{N}|{r'}_j|^{N_{j,k+1}}}{\Pi_{i=1}^{M}|r_i|^{M_{i,k+1}}}$. Thus, $l\in\mathbb N_*$  and  is  smaller than $m^*+n^*$.
Since   $\sum_{i=1}^M\overline{m_i}m_i-\sum_{j=1}^N\overline{n_j}n_j=1$, we have 
$\frac{\Pi_{i=1}^{M}|r_i|^{\overline{m_i}}}{\Pi_{j=1}^{N}|{r'}_j|^{\overline{n_j}}}=\gamma$ and  this says
$\frac{\Pi_{j=1}^{N}|{r'}_j|^{N_{j,k+1}}\Pi_{j=1}^{N}|{r'}_j|^{l\overline{n_j}}}{\Pi_{i=1}^{M}|r_i|^{M_{i,k+1}}\Pi_{i=1}^{M}|r_i|^{l\overline{m_i}}}=1$.
Thus
\begin{eqnarray*}
 (a+b)\Pi_{i=1}^{M}|r_i|^{M_{i,k}}&=&\Pi_{i=1}^{M}|r_i|^{-l\overline{m_i}}(a+b){\Pi_{i=1}^{M}|r_i|^{M_{i,k}-M_{i,k+1}}}\Pi_{i=1}^{M}|r_i|^{M_{i,k+1}}\Pi_{i=1}^{M}|r_i|^{l\overline{m_i}}\\&<&\Pi_{i=1}^{M}|r_i|^{-l\overline{m_i}}(a+b)\gamma^{-mn-m^*}\Pi_{i=1}^{M}|r_i|^{M_{i,k+1}}\Pi_{i=1}^{M}|r_i|^{l\overline{m_i}}. 
\end{eqnarray*}
Since $\Pi_{i=1}^{M}|r_i|^{M_{i,k}}=\gamma^{\sum_{i=1}^{M}m_iM_{i,k}}$ and $R< (a+b)\gamma^{\sum_{i=1}^{M}m_iM_{i,k}}$, we obtain
\begin{equation}\frac{\gamma^{mn+m^*+l\sum_{i=1}^M m_i\overline{m_i}}}{a+b}R< \Pi_{i=1}^{M}|r_i|^{M_{i,k+1}}\Pi_{i=1}^{M}|r_i|^{l\overline{m_i}}\end{equation}
Let  $T_1:= F_I\times F_J$ and $T_2:=(F_1^{l\overline{m_1}}\circ \dots \circ F_M^{l\overline{m_M}})\times ({F'_1}^{l\overline{n_1}}\circ \dots \circ {F'_N}^{l\overline{n_N}})$. Let $T=T_1\circ T_2$ and  $\alpha:=\Pi_{i=1}^{M}|r_i|^{M_{i,k+1}}\Pi_{i=1}^{M}|r_i|^{l\overline{m_i}}.$ Thus, 
$[DT]=\alpha\left[ \begin{smallmatrix} \pm1 & 0 
\\ 0 &\pm 1  
\end{smallmatrix} \right]$. 
\begin{itemize} 
\item If  $T(x,y)=\alpha(x,y)+(d,d')$, then  let $F(x)=\alpha x+e$, where $e=d-d'.$
To prove  $F(K-K')\subset K-K'$,  let $t\in K-K'$. There are $x\in K$ and $y\in K'$ such that 
$t=x-y$. Since,  $T(x,y)\in K\times K'$,  we have  $\alpha x+d\in K$ and $\alpha  y+d'\in K'$. Thus, $\alpha t +d-d'\in K-K'$.
Consequently, $F(t)\in K-K'$.
\item  If  $T(x,y)=-\alpha(x,y)+(d,d')$, then  let  $F(x)=-\alpha x+e$, where $e=d-d'.$  To prove  $F(K-K')\subset K-K'$,  let $t\in K-K'$. There are $x\in K$ and $y\in K'$ such that 
$t=x-y$. Since,  $T(x,y)\in K\times K'$,  we have  $-\alpha x+d\in K$ and $-\alpha  y+d'\in K'$. Thus, $-\alpha t +d-d'\in K-K'$.
Consequently, $F(t)\in K-K'$.
\item Otherwise, let $T_3=T\circ (F_K\times F_{K'})$. Thus 
$[DT_3]=\pm c \alpha\left[ \begin{smallmatrix} 1 & 0 
\\ 0 & 1  
\end{smallmatrix} \right]$.  Similar to cases above, there is the mapping  $F(x)=\pm c\alpha +e$  satisfying  $F(K-K')\subset K-K'$.
\end{itemize}
From $(2.1)$, we obtain
$$\frac{c\gamma^{mn+m^*+(m^*+n^*)\sum_{i=1}^M \overline{m_i}m_i}}{a+b}R|x-y|< |F(x)-F(y)|, ~~~~\forall x,y\in K-K'.$$
This completes the proof of the first assertion. 

For the  second assertion,  If  $ H^s(K-\lambda K')=0 $,  then  $  H^s(K-{r_i}^{-1}\lambda K')=0 $ and $  H^s(K-r'_j\lambda K')=0$. The first obtains since
  $$K-\lambda K'=\bigcup_{k=1}^{k=M}(r_kK+t_k)-\lambda K'=\bigcup_{k=1}^{k=M}\Big(\big(r_k(K-{r_k}^{-1}\lambda K')+t_k\big)\Big),$$
and  the second also has a similar proof.
On the other hand, affine  Cantor sets $K$ and $\frac{a}{b}K'$ have same size and the maps of the IFSs defining these sets are
orientation preserving.  By using   Theorem 2.12 in \cite{I} and  Theorem 2 in  \cite{PS}, we obtain   $H^s(K-\frac{a}{b}  K')=0$, since    $HD(K)+HD(K') \leq 1$  and $\frac{\log r_i}{ \log r_j'}\notin \mathbb{Q}$ for at least a pair $(i,j)$. Consequently,  $H^s(K- \frac{{r_j'}^{n}}{r_i^{m}} \frac{a}{b}  K' )=0$, for all $m,n\in \mathbb{N}$.  Thus,  the set of all  $\lambda$  that  $H^s(K-\lambda K')=0$    is  dense  in $\mathbb R^+$ since the set 
$\{  \frac{{r_j'}^{n}}{r_i^{m}} |~m,n\in \mathbb{N}\}$ is dense in $\mathbb R^+$.  The proof of Theorem \ref{A} is completed by replacing $K'$ with $b-K'$.
$~~\Box$
\section{Proof of Theorem  \ref{B} and its subsequent results} \label{sec3}
Suppose that $K$ is an affine Cantor set  with  Markov domain  $\{I_{0},I_{1},\cdots,I_{M}\}$
\label{sec2}  and  expanding map $\psi$ defined by $\psi_{|I_i}=:p_ix+e_i$ and also $K'$ is an affine Cantor set  with  Markov domain  $\{{I'}_{0},{I'}_{1},\cdots,{I'}_{N}\}$
\label{sec2}  and  expanding map $\psi'$ defined by $\psi'_{|I'_j}=:q_jx+f_j$. 
  Let $\Sigma:= \{0,1,2, \dots, M\},~\Sigma':=\{0,1,2, \dots, N\}$, and for   each $ i\in \Sigma$ and $j\in\Sigma'$, let  maps  $T_i$ and $T'_j$ from
$\mathbb{R}^* \times \mathbb{R}$ into itself given by
\begin{equation}
(s,t)\overset{T_i}\longmapsto
(p_is,p_it+e_i),\hspace{1cm}(s,t)\overset{T'_j}\longmapsto
(\frac{s}{q_j},t-\frac{f_j}{q_j}s).
\end{equation}
Each mapping of $(3.1)$ is called  a ``transferred renormalization operator"  or simply, an   `` operator"   \cite{HMP}. 
An important feature of these operators is that $T_i\circ T'_j=T'_j\circ T_i$.
 We say that  the  pair $(s,t)$ is a \textit{difference pair }for  $K$ and $K'$,  if there exists a relatively compact  sequence $(s_{n}, t_{n})$  in  the space  $\mathbb{R}^* \times
\mathbb{R}$ with $(s_1,t_1)=(s,t)$  such that the point  $(s_{n+1},t_{n+1})$ is obtained from $(s_{n}, t_{n})$ by  applying one of  operators (3.1). The following assertion  is proved in \cite{P1}:
\begin{itemize}
\item [$\bullet$] $t  \in K-sK'$   if and only if  $(s,t)$ is a difference pair for  $K$ and $K'$.
\end{itemize}
Let $\Pi_ \theta$ be the orthogonal projection from the plane to the line $l_\theta$  through the
origin making an angle $\theta$ with the positive $x$-axis.
For any non empty  set $E$  in the plane,  define the set $\cal{IP}$  as
$${\cal{IP}}={\cal{IP}}(E)=\{\theta\in[0,~\pi)| \pi_\theta \text{ is not one to one on $E$}  \}$$
In other words, $\cal{IP}$ is the set of all directions that can be obtained by connecting pairs of points in $E$. Eroglu  \cite{I}  proved that for two  affine Cantor sets  $K$ and  $K'$ in the real line with similarity dimensions $d$ and  $d'$, respectively, $H^{d+d'}(K+K')=0$, for Lebesgue almost all $\theta\in\cal{IP}$. If $d+d'<\frac{1}{2}$, then
$H^{d+d'}(K+K')>0$, for Lebesgue almost all $\theta\in \cal{IP}$. Recall that the similarity dimension of $K$  is the real number  $d$  such that $\sum_{i=1}^{M}|r_i|^d=1.$

\textbf{The proof of Theorem   \ref{B}:}
 Astels  extend Newhouse thickness for arbitrary compact sets \cite{A1}. He showed that for two non empty compact set $C_1$ and $C_2$,  if $\tau(C_1)\cdot\tau (C_2)<1$, then $$\tau(C_1+C_2)\geq\frac{\tau(C_1)+\tau(C_2)+2\tau(C_1)\tau(C_2)}{1-\tau(C_1)\tau(C_2)}.$$ If $\tau(C_1)\cdot\tau (C_2)\geq 1$,  the set  $C_1$  is  not contained in a gap of  $C_2$  and the set  $C_2$  is  not contained in a gap of  $C_1$, then $C_1\cap C_2=\emptyset$.
 We show that for two affine Cantor sets $K$ and $K'$ satisfying the condition $\tau(K)+\tau(K')+3\tau(K)\tau(K')\geq1$, the set $\cal{IP}$ contains an interval, and   generically ${\cal{IP}}=[0,~\pi)$. We can examine the option $\pi_ \theta(x,y)=x-\cot \theta y$ instead of $\Pi_ \theta(x,y)=x+\tan\theta y$, which is consistent with our notation.

\begin{lemma} \label{C}   Assume  that  $K$ and $ K'$  are two affine  Cantor sets such that $\tau(K)+\tau(K')+3\tau(K)\tau(K')\geq1$. Then 
\begin{itemize}
\item [$i)$]  the set  $\cal{IP}$ contains an interval  (indeed, $\tan^{-1} (\frac{b}{a})\in {\cal{IP}}^\circ$),
\item [$ii)$]   the set  $\cal{IP}$  generically (open dense with  full measure) is  the half interval $[0,~\pi)$.
\end{itemize}
\end{lemma}
\begin {proof}

For $(i)$, first consider the case where $p_0, q_0, p_M, q_N > 0$. Let  $C_1=[0,~\frac{a}{p_0}]\times  [0,~\frac{b}{q_0}]$ and  $C_2=[a-\frac{a}{p_M},~a]\times  [b-\frac{b}{q_N},~b]$, the  lower-left and upper-right corners  squares, respectively   of the first  step of construction of  $K\times K'$ in $[0,~a]\times [0,~b]$. For each $0< \theta< \frac{\pi}{2}$,   let difference sets,
$$K_1:=\pi_\theta(C_1\cap K\times K')=\frac{1}{p_0}K-\frac{\cot\theta}{q_0}K'$$ and $$K_2:=\pi_\theta(C_2\cap K\times K')=\frac{1}{p_M}K-\frac{\cot\theta}{q_N}K'+a-\frac{a}{p_M}-(b-\frac{b}{q_N})\cot\theta.$$

\begin{itemize} 
\item If $\tau(K)\tau(K')<1$, then  $K_1$ and $K_2$ have at least the  thicknesses $\frac{\tau(K)+\tau(K')+2\tau(K)\tau(K')}{1-\tau(K)\tau(K')}$.   The assumption of the theorem says $\tau(K_1)$ and $\tau(K_2)$ are greater than or equal to one.    For each $\phi_0:=\tan^{-1} (\frac{b-\frac{b}{q_N}}{a})\leq \theta\leq \theta_0:=\tan^{-1} (\frac{b}{a-\frac{a}{p_M}})$, the set  $\pi_\theta(C_1)\cap \pi_\theta(C_2)$ has non empty interior.  Newhouse gap lemma  gives $K_1\cap K_2\neq \emptyset$.
\item If $\tau(K)\tau(K')\geq 1$, then  there are the intervals $[\phi_0, ~\theta_0]$ containing  $\tan^{-1} (\frac{b}{a})$ and $(-t_0,~t_0)$ such that Cantor sets  $\frac{1}{p_0}K$  and  $\frac{\cot\theta}{q_0}K'+t$ does not fit into each other's   gaps,  for all $\theta\in [\phi_0, ~\theta_0]$ and $t\in(-t_0,~t_0)$. Gap lemma  gives $(-t_0,~t_0)\subset K_1$. Similarly,  there are the intervals $[\phi_0, ~\theta_0]$  containing  $\tan^{-1} (\frac{b}{a})$ and $(-t_0,~t_0)$  such that Cantor sets  $\frac{1}{p_M}K+a-\frac{a}{p_M}$  and  $\frac{\cot\theta}{q_N}K'+(b-\frac{b}{q_N})\cot\theta+t$ does not fit into each other's  gaps,  for all $\theta\in[\phi_0, ~\theta_0]$ and $t\in(-t_0,~t_0)$.  Note that 
$a=\max (\frac{1}{p_M}K+a-\frac{a}{p_M})=\max (\frac{\cot\theta}{q_N}K'+(b-\frac{b}{q_N})\cot\theta)$, when $\theta=\tan^{-1} (\frac{b}{a})$. Gap lemma  gives $(-t_0,~t_0)\subset K_2$. 
\end{itemize}
In both cases,  $K_1\cap K_2\neq \emptyset$, and this says $\pi_\theta$ is not one to one for each $\phi_0\leq \theta\leq \theta_0$. In general, let
$$K_1:=\psi^{-1}_{|I_0}(K)-\cot\theta{\psi'}^{-1}_{|I'_0}(K')  ~~~~~\text{and}~~~~~K_2:=\psi^{-1}_{|I_M}(K)-\cot\theta{\psi'}^{-1}_{|I'_N}(K'). $$
A proof similar to the one presented above can be given to show that $K_1 \cap K_2 \neq \emptyset$, and this completes the proof of assertion $(i)$.

Note that by substituting $b-K'$ for $K'$ in part $(i)$, we see that $\pi_\theta$ is not one to one for any $\theta$ in the interval $\phi_0 \leq |\theta| \leq \theta_0$.

For given   two  affine  Cantor sets $K$ and $K'$   with   convex hull  $[0,~a]$ and $[0,~b]$, respectively,  one can write $$K+ K'=a\big(\frac{1}{a}K+\frac{ b}{a} \frac{1}{b}K'\big),$$ in which $\frac{1}{a}K$ and $ \frac{1}{b}K'$  are affine  Cantor sets with convex hull  $[0,1]$.
Consequently, for given   affine Cantor sets $K_1$ and $K_2$ with arbitrary convex hulls,  there are  affine  Cantor sets $K$ and $K'$   with  convex hull  $[0,~1]$  and also  positive real  number $\lambda$ such that   $K_1+K_2$ is  an   affine translation of  $K+\lambda  K'$. Let,
$${\cal{K}}:=\big\{(K,s K')\big| ~~\text{   $K$ and $ K'$  have  convex hull  $[0,1]$}, ~~s>  0 \big\}$$
and
$${\cal{U}}:=\big\{(K,sK')\in{\cal{K}}\big| ~\exists~ m, n\in \mathbb{N}\cup\{0\}  ~\text{and} ~ {(i,j)\in\Sigma\times\Sigma'} \text{ s.t}~\cot \theta_0<|\frac{ p_i^m}{q_j^n}|s< \cot \phi_0\big\}.$$
 Obviously,  $\cal{U}$ contains  all  pairs $(K,s K')$ which $\frac{ \log| p_i|}{\log |q_j|}$  is an irrational number  for a $(i,j)\in\Sigma\times\Sigma'$.
It is not hard to show that $\cal{U}$ is an open and dense subset of $\cal{K}$.   Lastly, let $\cal{V}$ be the space of pairs $(K,K')$ of  affine   Cantor sets such that  there is   $(K,s K')\in{\cal{U}}$ which   $K_1+K_2$ is  an   affine translation of  $K+s  K'$. Clearly, $\cal{V}$  is an open and dense  subset  in  the space     of  all pairs of   affine Cantor sets. 

For $(ii)$,  suppose that $(K,s_0K')\in {\cal{U}}$ with $s_0>0$ are chosen.  Thus, there are $ m, n\in \mathbb{N}\cup\{0\} $ and $ (i,j)\in\Sigma\times\Sigma'$  such that $\cot \theta_0<|\frac{ p_i^m}{q_j^n}|s_0< \cot \phi_0$.  From above,  there are two disjoint points $(x_1,y_1),(x_2,y_2)\in K\times K'$  such that $\pi_\theta(x_1,y_1)=\pi_\theta(x_2,y_2),$   where   $\cot \theta=\frac{ p_i^m}{q_j^n}s_0$ . Let  $s=\frac{ p_i^m}{q_j^n}s_0$ and $t=\pi_\theta(x_1,y_1)=x_1-sy_1$,  thus, the  pair $(s,t)$ is a difference pair for  $K$ and $K'$. This gives  the existence of two  disjoint  relatively compact  sequence $(s_{n}, t_{n})$ and $(s'_{n}, t'_{n})$  in  the space  $\mathbb{R}^* \times
\mathbb{R}$ with $(s_1,t_1)=(s,t)=(s'_{1}, t'_{1})$  such that the point  $(s_{n+1},t_{n+1})$ and $(s'_{n+1},t'_{n+1})$ are  obtained from $(s_{n}, t_{n})$ and  $(s'_{n},t'_{n})$, respectively, by  applying one of  operators (3.1). On the other hand, there is  an interval $I\subset [-s_0,~1]$, if $s_0>0$ or $I\subset [0,~1+s_0]$, if $s_0<0$, such that
$$T_i^n\circ{T'_j}^m\big(\{s_0\}\times I\big)=\{\frac{ p_i^m}{q_j^n}s_0\}\times[-\frac{ p_i^m}{q_j^n}s_0,~1],$$ or
$$T_i^n\circ{T'_j}^m\big(\{s_0\}\times I\big)=\{\frac{ p_i^m}{q_j^n}s_0\}\times[0.~1+\frac{ p_i^m}{q_j^n}s_0],$$ respectively. Thus, there exists a $t_0\in I$, such that $T_i^n\circ{T'_j}^m(s_0,t_0)=(s,t)$.
Thus, there are  two  disjoint  relatively compact  sequence $(\tilde s_{n}, \tilde t_{n})$ and $( \tilde s'_{n},\tilde t'_{n})$  in  the space  $\mathbb{R}^* \times
\mathbb{R}$ with $(\tilde s_1,\tilde t_1)=(s_0,t_0)=(\tilde s'_{1}, \tilde t'_{1})$  such that the point  $(\tilde s_{n+1},\tilde t_{n+1})$ and $(\tilde s'_{n+1},\tilde t'_{n+1})$ are  obtained from $(\tilde s_{n},\tilde t_{n})$ and  $(\tilde s'_{n},\tilde t'_{n})$, receptively by  applying one of  operators (3.1).  These two sequence gives  the existence  two disjoint points $(\tilde x_1,\tilde y_1),(\tilde x_2,\tilde y_2)\in K\times K'$ such that $t_0=\tilde x_1-s_0\tilde y_1=\tilde x_2-s_0\tilde y_2$. This says that the projection map   $\pi_ \theta$  is not one to one on the $ \cal{U}$.  Thus $\Pi_ \theta$   is not one to one on the $\cal{V}$  and this completes the proof of the lemma.
\end{proof}

Now, suppose that  for at least one pair $(i, j)$ we have $\frac{\log |r_i|}{ \log |r_j'|}\notin \mathbb{Q}$.  Then the relation $s:=HD\big(K-\cot \theta K'\big)= d+d'$ holds for all $\theta\in (0,~\pi)-\frac{\pi}{2}$ \cite{PS}. On the other hand, for almost all ${\theta \in \mathcal{IP}}(K\times K')$, we have $H^{d+d'}(K+K')=0$ (in the sense of Lebesgue measure) \cite{I}. From Lemma \ref{C}, the set ${\cal{IP}}(K\times K')$ has full  Lebesgue measure in the interval $[0, \pi)$, since $\tau(K)+\tau(K')+3\tau(K)\tau(K')\geq1$. Therefore, $H^s(K-\lambda K')=0$ holds for almost all $\lambda \in \mathbb R$, and the proof of Theorem \ref{B} is complete.   $\Box$

Recall of  \cite{P1}  that a homogeneous  Cantor set is middle, if  
 the removed intervals at the first step of its construction process  have equal lengths. Suppose that the underlying  homogeneous  Cantor set $K$ (resp. $K'$) is  middle    with  convex hull  $[0,~a]$ (resp. $[0,~b]$), and the length of  the removed intervals  in the first step from its construction is  $c$ (resp. $c'$). We showed in the case   $\tau(K)\cdot\tau(K')<1$.
  \begin{itemize}
  	\item [I)]   if $\frac{\log \lambda}{\log \mu}$ is  irrational, then  $K+ K'\neq [0,~a+b]$,
  	
  	\item [II)] if $\frac{\log \lambda}{\log \mu}=\frac{n}{m}$  with $(m,n)=1$ and $\gamma:=\lambda^{-\frac{1}{n}}$,  then 
  	\begin{itemize}
  		\item [$\bullet$]  the condition  $\tau(K)\cdot\tau(K')<\frac{1}{\gamma}$ gives   $K+ K'\neq [0,~a+b],$ 
  		
  		\item [$\bullet$]  the condition $\tau(K)\cdot\tau(K')\geq \frac{1}{\gamma}$ gives $K+ K'=[0,~a+b]$  if and only if
  		\begin{align*}
  		\big\lceil n\log_\lambda\frac{\lambda a}{c'}\big\rceil=\big\lfloor n\log_\lambda\frac{c}{\mu b}\big\rfloor+1,
  		\end{align*} 
  	\end{itemize} 
  \end{itemize} 
 where $\lfloor ~\rfloor$ and $\lceil~\rceil$ are floor and ceiling functions, respectively. 
 Take  natural number  $M$ and positive real number $\lambda$ with $\lambda <\frac{1}{M+1}$. The \textit{$M$-middle}   Cantor set $C_\lambda^M$ is defined as  a middle  homogeneous  Cantor set with  convex hull  $[0,~1]$ such that  the number of  removed intervals in the first step of its construction process  is $M$ and  they have the same  length  $\frac{1-(M+1)\lambda}{M}$. Algebraically,
 $$C_\lambda^M:=\Big\{\frac{1-\lambda}{M}\sum_{n=0}^{\infty} a_n\lambda^n\big|~~a_n\in \{0,1,\cdots, M\}\Big\}.$$
 When $M=1$, the set $C_\lambda^1$ is the same   middle-$(1-2\lambda)$  Cantor set  $C_\lambda$.
Obviously, $HD(C_\lambda^M)=-\log_\lambda (M+1)$ and $\tau(C_\lambda^M)=\frac{M\lambda}{1-(M+1)\lambda}$. 
Among our results, we also showed that if $M=N$ and $m\neq n$,  then there is  $\lambda_0 <\lambda_1$, where 
$\lambda_1$ satisfying   $\tau(C_{\lambda_1}^M)\cdot\tau(C_{\lambda_1^\frac{m}{n}}^N)=1$, such that for all  $\lambda_0 \leq\lambda\leq \lambda_1$ we  always have  $C_\lambda^2+C_{\lambda^\frac{m}{n}}^2=[0,~2]$. \\
Let  $C_\lambda$  and $C_\mu$  be the middle-$(1-2\lambda)$ and  middle-$(1-2\mu)$ Cantor
sets, respectively. The relation  $\tau(C_\lambda)+\tau(C_\mu)+3\tau(C_\lambda)\tau(C_\mu)\geq1$  is equivalent to  $\frac{2\lambda}{1-3\lambda}\frac{2\mu}{1-3\mu}=\tau(C^2_\lambda)\cdot\tau(C^2_\mu)\geq 1$. From Lemma   \ref{C}, the set  $\cal{IP}=:\cal{IP}_{(\lambda,\mu)}$ contains an interval around $\frac{\pi}{4}$  if $\tau(C^2_\lambda)\cdot\tau(C^2_\mu)\geq 1$.

\begin{proposition} \label{D} For all $m,n\in\mathbb N$ with $m\neq n$,  there are  smooth curves
\begin{align*}
r_{m,n}: ~I_{m,n}\longmapsto R\\
\gamma \rightarrow \big(\gamma^{-n},\gamma^{-m} \big)
\end{align*}
where $R:=\Big\{(\lambda,\mu)\Big|~\tau(C^2_\lambda)\cdot\tau(C^2_\mu)<1\Big\},$
 the   ${\cal{IP}}_{r_{m,n}(\gamma)}$ contains
an interval  around $\frac{\pi}{4}$.
\end{proposition}
\begin{proof}
Take the notation of the proof of  Theorem   \ref{B}. For each $0< \theta< \frac{\pi}{2}$, we have  $K_1=\lambda C_\lambda-\mu\cot \theta  C_\mu$ and 
 $K_2=\lambda C_\lambda-\mu\cot \theta  C_\mu+1-\lambda- (1-\mu)\cot \theta$. In these cases,  $K_2$  is an  affine copy of $K_1$ and this   gives
 \begin{eqnarray*}
K_2-K_1&=&\lambda( C_\lambda- C_\lambda)- \mu\cot \theta(C_\mu-C_\mu)+1-\lambda- (1-\mu)\cot \theta\\
&=&\lambda( C_\lambda+ C_\lambda)- \mu\cot \theta(C_\mu+C_\mu)+1-2\lambda- (1-2\mu)\cot \theta.
\end{eqnarray*}
 The set  $\lambda( C_\lambda+ C_\lambda)$ is a   middle   homogeneous Cantor set with convex hull $[0, 2\lambda]$ and a Markov partition consisting of three intervals of length $2\lambda^2$.
 The slopes of the expanding maps are equal to $\lambda$ and so  $c_\lambda=\lambda(1-3\lambda)$.  Indeed,  $\lambda( C_\lambda+ C_\lambda)=2\lambda C^2_\lambda$, also, $\tau\big( C^2_\lambda\big)=\frac{2\lambda}{1-3\lambda}$.
Similarity,  $ \mu\cot \theta(C_\mu+C_\mu)=2 \mu\cot \theta C^2_\mu$ is a   middle    homogeneous Cantor set with  convex hull  $[0,~2\mu\cot \theta]$ and the slopes of the expanding maps are equal to $\mu$  and  the  length of  gaps are $c_\mu=\mu(1-3\mu)\cot \theta$. Also, $\tau\big( C^2_\mu\big)=\frac{2\mu}{1-3\mu}$.

If $\frac{\log \lambda}{\log \mu}=\frac{n}{m}$  with $(m,n)=1$ and $\gamma:=\lambda^{-\frac{1}{n}}$,  from above  $ \frac{1}{\gamma}\leq\frac{2\lambda}{1-3\lambda}\frac{2\mu}{1-3\mu}<1$ gives 
\begin{eqnarray*}
K_2-K_1&=&[-2 \mu\cot \theta,~2\lambda]+1-2\lambda- (1-2\mu)\cot \theta\\
&=&[1-2\lambda-\cot \theta,~1- (1-2\mu)\cot \theta],
\end{eqnarray*}
 if and only if
$$\big\lceil n\log_\lambda\frac{2\lambda^2 }{\mu(1-3\mu)\cot \theta}\big\rceil=\big\lfloor n\log_\lambda\frac{\lambda(1-3\lambda)}{2\mu^2\cot \theta}\big\rfloor+1.$$ 
or equivalently, 
\begin{align}\big\lceil n\log_\lambda\frac{2\lambda }{(1-3\mu)\cot \theta}\big\rceil=\big\lfloor n\log_\lambda\frac{(1-3\lambda)}{2\mu\cot \theta}\big\rfloor+1.\end{align}
On the other hand,  since  $m\neq n$ and $M=N=2$,  there is  $\lambda_0 <\lambda_1$, where 
$\lambda_1$ satisfying  $\frac{2\lambda_1}{1-3\lambda_1}\frac{2\lambda_1^\frac{m}{n}}{1-3\lambda_1^\frac{m}{n}}=1$  such that for all  $\lambda_0 \leq\lambda\leq \lambda_1$ we  always have  $C_\lambda^2+C_{\lambda^\frac{m}{n}}^2=[0,~2]$.  Again from above,    $C_\lambda^2+ C_\mu^{2}=[0,~2],$  if and only if
$$\big\lceil n\log_{\lambda}\frac{2\lambda}{1-3\lambda^{\frac{m}{n}}}\big\rceil=\big\lfloor n\log_{\lambda}\frac{1-3\lambda}{2\lambda^{\frac{m}{n}}}\big\rfloor+1.$$
 Among the proofs of the paper \cite{P1}, we see that   $n\log_{\lambda_1}\frac{2\lambda_1}{1-3\lambda_1^{\frac{m}{n}}}\notin\mathbb{Z}$, since $\lambda_0\neq \lambda_1$. This says $ n\log_{\lambda_1}\frac{1-3\lambda_1}{2\lambda_1^{\frac{m}{n}}}\notin\mathbb{Z}$, since   $\frac{2\lambda_1}{1-3\lambda_1}\frac{2\lambda_1^\frac{m}{n}}{1-3\lambda_1^\frac{m}{n}}=1$.  Thus for each $\gamma\in I_{m,n}:=[{\lambda_0}^{\frac{-1}{n}},~{\lambda_1}^{\frac{-1}{n}}]$, there is an interval   $(\frac{\pi}{4}- \epsilon,~\frac{\pi}{4}- \epsilon)$ such that the equality $(3.2)$ works for  the pair $I_{m,n}(\gamma)=(\lambda,\mu)$, for all $ \theta \in (\frac{\pi}{4}- \epsilon,~\frac{\pi}{4}- \epsilon)$. Consequently, $K_2-K_1$ is the same closed interval mentioned above.  
Since $\pi_{\frac{\pi}{4}}(K_2\times K_1)=[-2\lambda,~2\mu]$, there is  $0<\epsilon_1\leq \epsilon$ such that $0\in K_2-K_1$, for all $\theta \in (\frac{\pi}{4}- \epsilon_1,~\frac{\pi}{4}- \epsilon_1)$.
This  gives $K_1\cap K_2\neq \emptyset$, for all $\theta \in (\frac{\pi}{4}- \epsilon_1,~\frac{\pi}{4}- \epsilon_1)$. This says $\pi_\theta$ is not one to one for each $\theta \in (\frac{\pi}{4}- \epsilon_1,~\frac{\pi}{4}- \epsilon_1)$ and  so  the set  ${\cal{IP}}_{r_{m,n}(\gamma)}$ contains an interval  around $\frac{\pi}{4}$  for all  $\lambda_0 \leq\lambda\leq \lambda_1$. This completes the proof of the proposition.
\end{proof}
At end, we   study topological  and measure  theoretic structures of  $C_\lambda+C_\mu$ for  all $(\lambda,\mu)$
  in the region $[0,~\frac{1}{2}]\times[0,~\frac{1}{2}]$  in the form that Eroglu studied  in  \cite{I}. 
For $(\lambda,\mu)$ has been indicated as region $\text{ III}:=\{(\lambda,\mu)|~ \tau(C_\lambda)\cdot\tau(C_\mu)\geq1\}$ in  Figure 1, we have  $C_\lambda+C_\mu=[0,~2]$.
The curve $d_\lambda+d_\mu=1$, where $d_\lambda=HD(C_\lambda)$ and  $d_\mu=HD(C_\mu)$,   represents the lower boundary of region   $\text{ II}:=\{(\lambda,\mu)|~d_\lambda+d_\mu\geq1,~~ \tau(C_\lambda)\cdot\tau(C_\mu)<1\}$ in Figure 1. In this region, for almost all parameter pairs on each vertical or horizontal line we have
    $|C_\lambda+C_\mu|>0$.
In region  $\text{I}$ which  divide it into subregions 
\begin{itemize}
  		\item [$\bullet$]  $ \text{Ia}:=\{(\lambda,\mu)|~d_\lambda+d_\mu<1,~~~~ \tau(C^2_\lambda)\cdot\tau(C^2_\mu)\geq 1\},$
  		\item [$\bullet$]  $\text{Ib}:=\{(\lambda,\mu)|~\tau(C^2_\lambda)\cdot\tau(C^2_\mu)< 1,~~~~d_\lambda+d_\mu\geq\frac{1}{2}\},$
  		\item [$\bullet$]  $ \text{Ic}:=\{(\lambda,\mu)|~d_\lambda+d_\mu<\frac{1}{2}\},$
  	\end{itemize} 
\begin{figure}[ht]
 	\centering 
 	\scalebox{0.7} 
 	{\includegraphics {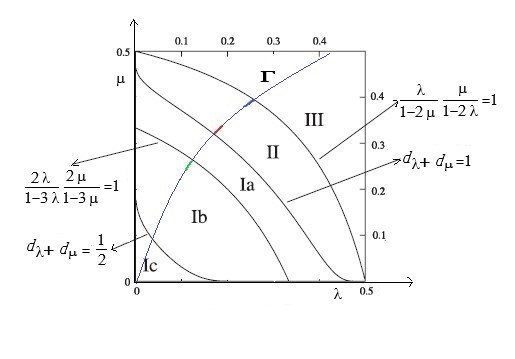}}\caption{\small{Different cases for $(C_\lambda,C_\mu)$ when the two sets are middle  Cantor sets.}}
 	\label{fig:exm} 
 \end{figure} 
we have  $d_\lambda+d_\mu<1$, and so   $|C_\lambda+C_\mu|=0$. In region  $\text{Ia}$,    the set  $\cal{IP}_{(\lambda,\mu)}$ contains an interval. In region $\text{Ic}$, Lebesgue almost all projections have positive ($d_\lambda+d_\mu$ )-dimensional Hausdorff measure.  For all $(\lambda, \mu)$ on the line $\Gamma: \mu=\lambda$ located in the region $\text{I}\cup\text{II}$ (indeed, $\lambda<\frac{1}{3}$), we have $|C_\lambda+C_\mu|=0$. For $(\lambda, \mu)$ with $\lambda > 1/6$ (in which a portion of the line lies in region $\text{Ib}$), the set $\cal{IP}_{(\lambda,\mu)}$ contains  an interval, and if $\lambda < 1/9$, then ${HD\big(\cal{IP}_{(\lambda,\mu)}\big)} < 1$.

Take  different  natural numbers $m$ and $ n$ with $(m,n)=1$  and let the curve $\Gamma:~\mu=\lambda^\frac{m}{n}$. The portion of the curve $\Gamma$ that lies in region II and connects to the boundary of region III (indicated in blue in Figure 1) satisfies the condition $C_\lambda + C_\mu = [0, ~2]$.  The portion of the curve $\Gamma$ that lies in region $\text{II}$ and connects to the boundary of region $\text{I}$ (marked in red in Figure 1) satisfies the condition $|C_\lambda + C_\mu| = 0$ (references can be found in   \cite{I}). The result we add here (Proposition \ref{D}) is that there exists a portion of the curve $\Gamma$ lying in $\text{Ib}$ connected to the boundary of  region  $\text{Ia}$ and (shown in green in the figure) along which the set $\cal{IP}_{(\lambda,\mu)}$ contains an interval   around $\frac{\pi}{4}$.


\bigskip
\begin{thebibliography}{99}

\bibitem{A} R. Anisca, C. Chlebovec, M. Ilie, \textit{The structure of arithmetic sums of affine Cantor sets}, Real Anal. Exch.   (37)2 (2012), 325-332.
\bibitem{A1}  S. Astels,   \textit{Cantor sets and numbers with restricted partial quotients},  Trans. Amer. Math. Soc.
(352)1 (1999),   133-170.
\bibitem{Fa} K.  Falconer, \textit{Fractal geometry: mathematical foundations and applications}, Wiley, (2003).
\bibitem{Fa1} K.  Falconer, \textit{Techniques in fractal geometry}, Wiley,  (1997).
\bibitem{HMP} B. Honary, C. G. Moreira,  M. Pourbarat, \textit{Stable
intersections of affine Cantor sets}, Bull. Braz. Math. Soc. (36)3 (2005),  363-378.
\bibitem{I} K. Ilgar Eroglu, \textit{On the arithmetic sums of Cantor sets}, Nonlinearity (20)5 (2007),  1145-1161.
\bibitem{Ma} A. McDonald, K. Taylor, \textit{Finite point configurations in products of thick Cantor sets and a robust nonlinear Newhouse Gap Lemma}, Math. Proc. Cambridge Philos. Soc. (175)2 (2023),  285 - 301.
\bibitem{MO} P. Mendes,  F. Oliveira,  \textit{On the topological structure of the arithmetic sum of two Cantor sets},
 Nonlinearity (7)2 (1994),  329-343.
\bibitem{M} C. G. Moreira, \textit{Stable intersections of Cantor sets and homoclinic bifurcations}, Ann. Inst. H. Poincare Anal. Non Lineaire (13)6 (1996), 741-781.
\bibitem{M2} C. G. Moreira, \textit{Geometric properties of images of cartesian products of regular Cantor sets by differentiable real maps}, Math. Z. (303)3 (2023).
\bibitem{M3} C. G. Moreira, E. M.  Munoz, \textit{Sums of Cantor sets whose sum of dimensions is close to 1}, Nonlinearity,  (16)5 (2003), 1641–1647.
\bibitem{M4} C. G. Moreira, E. M.  Munoz, J.  Rivera-Letelier \textit{On the topology of arithmetic sums of regular Cantor
sets }, Nonlinearity,  (13)6 (2000), 2077–2087.
\bibitem{MY}C. G. Moreira, J.-C. Yoccoz, \textit{Stable intersections of regular Cantor
sets with large Hausdorff dimension}, Ann. of Math. (154)1 (2001), 45-96.
\bibitem{P} M. Pourbarat, \textit{Stable intersection of middle-$\alpha$
Cantor sets}, Comm. in Cont. Math. (17)5 (2015).
\bibitem{P3} M. Pourbarat, \textit{Stable intersection  of  affine  Cantor sets defined by  two expanding maps}, J. Differential Equations  (266)4 (2019),  2125-2141.
\bibitem{P1}  M. Pourbarat, \textit{Topological structure of  the sum  of  two  homogeneous  Cantor sets},  Ergod. Theory Dyn. Syst.  (43)5 (2023), 1712 - 1736.
\bibitem{P11} M. Pourbarat, \textit{On the arithmetic difference   of
middle  Cantor sets}, Discrete Contin. Dyn. Syst,  38 (2018). 
\bibitem{P4} M. Pourbarat, \textit{On the sum   of two
homogeneous  Cantor sets},   Discrete Contin. Dyn. Syst.  (45)6 (2025),  1745-1766.
\bibitem{PS} Y. Peres,  P. Shmerkin, \textit{Resonance between Cantor
sets}, Ergod. Theory Dyn. Syst.  (29)1  (2009), 201-221.
\bibitem{T1}  Y. Takahashi,   \textit{Sums of two homogeneous Cantor sets},  Trans. Amer. Math. Soc.
(372)3 (2019),  1817–1832.
\bibitem{PT}  J. Palis,  F. Takens, \textit{Hyperbolicity and sensitive chaotic dynamics at
homoclinic bifurcations}, Cambridge Univ. Press, Cambridge, (1993).
\bibitem{PY} J. Palis,  J.-C. Yoccoz, \textit{On the arithmetic sum of
regular Cantor sets}, Ann. Inst. Henri Poincare 14 (1997),  439-456.
\bibitem{S1} A. Sannami, \textit{An example of a regular Cantor set whose
difference set is a Cantor set with positive measure}, Hokkaido
Math. J.  21 (1992), 7-24.
\bibitem{S2} B. Solomyak, \textit{On the arithmetic sums of  Cantor sets},
Indag. Mathem (1997), 133-141.
\bibitem{Y} A. Yavicoli,  \textit{Patterns in thick compact sets}, Israel J. Math. (244)1 (2021), 95-126.
\end{thebibliography}
\end{document}